\documentclass[12pt, reqno]{amsart}

\usepackage{mathtools}
\mathtoolsset{showonlyrefs}
\numberwithin{equation}{section}

\usepackage[dvipdfmx]{graphicx}

\usepackage{mathtools}
\mathtoolsset{showonlyrefs}

\usepackage{amssymb}
\usepackage{amsthm}
\usepackage{amscd}
\usepackage{amsmath}
\usepackage{epic,eepic}
\usepackage{mathrsfs}
\usepackage{latexsym}
\usepackage{pifont}
\usepackage[all]{xy}
\usepackage {cancel}

\theoremstyle{plain}
\newtheorem{lemma}{Lemma}[section]
\newtheorem{prop}[lemma]{Proposition}
\newtheorem{thm}[lemma]{Theorem}
\newtheorem{cor}[lemma]{Corollary}
\newtheorem{intthm}{Theorem}

\theoremstyle{definition}

\newtheorem{rem}[lemma]{Remark}
\newtheorem{defi}[lemma]{Definition}
\newtheorem{exa}[lemma]{Example}

\newtheorem{problem}{Problem}

\newcommand{\bde}{\begin{defi}}
\newcommand{\ede}{\end{defi}\vspace{1mm}}
\newcommand{\ble}{\begin{lemma}}
\newcommand{\ele}{\end{lemma}}
\newcommand{\bpr}{\begin{prop}}
\newcommand{\epr}{\end{prop}}
\newcommand{\bt}{\begin{thm}}
\newcommand{\et}{\end{thm}}
\newcommand{\bco}{\begin{cor}}
\newcommand{\eco}{\end{cor}}
\newcommand{\bre}{\begin{rem}}
\newcommand{\ere}{\end{rem}}
\newcommand{\bex}{\begin{exa}}
\newcommand{\eex}{\end{exa}}
\newcommand{\bpf}{\begin{proof}}
\newcommand{\epf}{\end{proof}}

\newcommand{\mcD}{\mathcal{D}}
\newcommand{\mcE}{\mathcal{E}}
\newcommand{\mcF}{\mathcal{F}}
\newcommand{\mcG}{\mathcal{G}}
\newcommand{\mcH}{\mathcal{H}}

\newcommand{\mcK}{\mathcal{K}}
\newcommand{\mcL}{\mathcal{L}}
\newcommand{\mcM}{\mathcal{M}}
\newcommand{\mcN}{\mathcal{N}}
\newcommand{\mcO}{\mathcal{O}}

\newcommand{\mcS}{\mathcal{S}}
\newcommand{\mcT}{\mathcal{T}}

\newcommand{\mbC}{\mathbb{C}}

\newcommand{\mbH}{\mathbb{H}}

\newcommand{\mbP}{\mathbb{P}}

\newcommand{\mbZ}{\mathbb{Z}}

\newcommand{\SSP}{\vspace{3mm}}
\newcommand{\LSP}{\vspace{5mm}}
\newcommand{\mr}{\mathrm}
\newcommand{\M}{m}
\newcommand{\V}{\mcE}
\newcommand{\N}{N}

\title[Complete flags on flat vector bundles  in positive characteristic]{Complete flags on  flat  vector bundles  \\    in positive characteristic}
\author{Youhei Morita}
\author{Yasuhiro Wakabayashi}

\address{\emph{Youhei Morita}
 \newline
 \textnormal{Graduate School of Information Science and Technology, Osaka University, Suita, Osaka 565-0871, JAPAN.}
 \newline
 \textnormal{\texttt{morita.youhei@ist.osaka-u.ac.jp}}}

 \address{\emph{Yasuhiro Wakabayashi}
 \newline
 \textnormal{Graduate School of Information Science and Technology, Osaka University, Suita, Osaka 565-0871, JAPAN.}
 \newline
  \textnormal{\texttt{wakabayashi@ist.osaka-u.ac.jp}}}

\begin{document}
\maketitle


\begin{abstract}
Let $X$ be a connected, smooth, and projective curve of genus $g$ over an algebraically closed field of characteristic $p>0$. This paper investigates a characteristic-$p$ analogue of a well-known fact concerning flat vector bundles in characteristic zero. That is to say, we prove that the inequality $g\le1$ holds if and only if any flat vector bundle on $X$ admits a complete flag. We also explore a generalization of this result in a broader setting using Berthelot's differencial operators of level $m$.
\end{abstract}
\tableofcontents

\section{Introduction} \label{S1}
\LSP

Let $X$ be a connected, smooth, and projective curve of genus $g\ge0$ over an algebraically closed field $k$. By a flat vector bundle on $X$, we mean a pair $(\V,\nabla)$, consisting of a vector bundle $\V$ on $X$ and a (flat) connection $\nabla$ on it.
When the base field $k$ is the field of complex numbers $\mbC$,
a flat vector bundle can be regarded as a $\mcD$-module,
where $\mcD$ is the sheaf of differential operators on $X$, and corresponds,
via the most basic form of the Riemann-Hilbert correspondence, to a group representation 
\begin{align}
\rho_{(\V,\nabla)}:\pi_{1}(X^\mr{an}) \rightarrow \mr{GL}_r (\mbC)
\end{align}
uniquely given up to conjugation,
where $r:=\mathrm{rank}(\V)$ and $\pi_{1}(X^{an})$ denotes the fundamental group of the associated Riemann surface $X^\mr{an}$ (with respect to a fixed base point).
Flat vector bundles and their associated group representations have been extensively studied in various contexts, including non-Abelian Hodge theory and the geometric Langlands program.

Now, let us focus on a curve of genus $g=1$; see ~\cite{Bis2}, ~\cite{FaLo}, ~\cite{FLM},  ~\cite{Kaw}, and (a part of) ~\cite{Wak8} for studies on such objects.
One of the fundamental  facts about these bundles when $k = \mbC$
  is the existence of a {\it complete flag} (i.e., a {\it splitting}, in the terminology of  ~\cite{Ati1},  ~\cite{Ati2}) preserved by  a given connection.
Specifically, for  any flat vector bundle $(\V, \nabla)$ on the complex curve $X/\mbC$ with $g =1$,
there exists  a sequence of subbundles
\begin{align} \label{eq99}
0 = \V^0 \subseteq  \V^1 \subseteq \cdots \subseteq \V^r = \V
\end{align}
 such that, for every $j=1, \cdots, r$, the subquotient $\V^j /\V^{j-1}$ is a line bundle and $\V^j$ is closed under $\nabla$.
In particular, 
$(\V, \nabla)$ can be obtained as a successive extension of flat line bundles, which  allows for a detailed understanding of  the moduli space of flat vector bundles on $X$.

The proof of this fact is elementary.
Indeed, 
since $\pi_1 (X^\mr{an}) \cong \mbZ^2$, the image of the group representation $\rho_{(\V, \nabla)}$ can be simultaneously upper triangularizable,
which implies  the desired result.
The same statement holds  for  genus-$0$ curves, as  their  fundamental groups are trivial.
On the other hand, when the underlying curve is hyperbolic, i.e., $g >1$,  there exist many irreducible group  representations $\pi_1 (X^\mr{an}) \rightarrow \mr{GL}_r (\mbC)$, and  the corresponding flat vector bundles do not admit complete flags  as in \eqref{eq99}.
It follows that {\it the existence of such flags depends entirely on the genus of the underlying curve}, and a standard base-change argument shows that this fact generalizes to the case where $k$ is an arbitrary algebraically closed field of characteristic zero.

Then, a natural question is whether a similar result holds in characteristic $p > 0$.
To the best of our knowledge, the analogous result  in positive characteristic has not been established, and  the proof outlined above does not apply in this setting because   it fundamentally relies on the analytic nature of the Riemann-Hilbert correspondence.
Moreover, some flat vector bundles in characteristic $p$ cannot be obtained by reducing those in characteristic zero modulo $p$.
Therefore, to establish a positive characteristic analogue  of the aforementioned result,
we need to address  phenomena that  are specific to characteristic $p$.

Let us suppose that $k$ has characteristic $p$.
We should note that 
there are different  variations   of 
the sheaf of differential operators ``$\mcD$" on $X$.
 On the one hand,  the sheaf $\mcD^{(\infty)}$ of differential operators  in the sense of A. Grothendieck was discussed  in 
 ~\cite{EGA4}.
 A $\mcD^{(\infty)}$-module is often referred to as a {\it stratified sheaf} and can be interpreted as an $\mcO_X$-module admitting infinite Frobenius descent.
On the other hand,  the sheaves  $\mcD^{(\M)}$ of differential operators of {\it level $\M \in \mbZ_{\geq 0}$}  were introduced by P. Berthelot (cf.  ~\cite{PBer1}, ~\cite{PBer2}) as  essential ingredients in defining arithmetic $\mcD$-modules.
A $\mcD^{(0)}$-module structure  amounts essentially to a choice of  a flat connection in the usual sense.

As described below, our  main result provides  a necessary and sufficient condition for   the existence of a complete flag on a $\mcD^{(\M)}$-bundle (i.e., a vector bundle equipped with a $\mcD^{(\M)}$-action).

\SSP
\begin{intthm}[cf. Theorems \ref{Th4},  \ref{Th2}, and \ref{Cor32}] \label{ThA}
Let $\M$ be an element of $\mbZ_{\geq 0} \sqcup \{ \infty \}$ and  $r$  a positive integer.
Then, 
the following two conditions (a) and (b) are equivalent:
\begin{itemize}
\item[(a)]
The inequality $g \leq 1$ holds;
\item[(b)]
 For any given rank $r$, any $\mcD^{(\M)}$-bundle $(\V, \nabla)$ of rank $r$ admits a complete flag, i.e., there exists  a sequence of subbundles  
 \begin{align} \label{Eq13}
 0 = \V^0 \subset \V^1 \subseteq  \cdots \subseteq \V^r = \V
 \end{align} 
 of $\V$ such that, for every $j=1, \cdots, r$,   the subquotient $\V^j/\V^{j-1}$ is a line bundle and $\V^j$ is preserved by $\nabla$.
 \end{itemize}
  \end{intthm}
  
  \begin{rem}
As pointed out by an anonymous referee, it is highly expected that a more general reduction theory holds for principal $G$-bundles, where $G$ is a geometrically connected, reductive algebraic group over an algebraically closed (or perfect) field such that $G$ is of low height with respect to the characteristic $p > 0$. While we focus on vector bundles (i.e., $\mr{GL}_r$-bundles) in this paper, developing the precise analogue for such principal $G$-bundles is an important direction for future research.
\end{rem}


\LSP
\section{Preliminaries} \label{S2}
\LSP

In this subsection,
we recall some fundamental aspects of connections and $\mcD$-modules in positive characteristic.
For previous studies on these objects, we refer the reader to, e.g.,  ~\cite{PBer1}, ~\cite{PBer2}, \, 
~\cite{Kat}, and ~\cite{Wak8}.
 
\subsection{Notation and Conventions} \label{SS7}
 
 Throughout this paper, we fix a prime number $p$ and  an algebraically closed field $k$ of characteristic $p$.
We also  fix 
a connected, smooth, and projective curve $X$ over $k$ of genus $g \in \mbZ_{\geq 0}$.
Denote by $\Omega_X$ (resp., $\mcT_X$) the cotangent bundle  (resp., the tangent bundle) of $X/k$ and by $d$ the universal derivation $\mcO_X \rightarrow \Omega_X$.

Denote by $F_X$ the absolute Frobenius endomorphism of $X$.
 For a positive  integer $\M$,
 the {\bf $\M$-th Frobenius twist} of $X$ over $k$ is, by definition,
 the base-change $X^{(\M)}$ of $X$  along the $\M$-th iterate of the absolute Frobenius automorphism of $\mr{Spec}(k)$.
The morphism $F_{X/k}^{(\M)} : X \rightarrow X^{(\M)}$ induced naturally from $F_X^{\M }$ is called the {\bf $\M$-th relative Frobenius morphism} of $X$ over $k$.
For convenience, 
we write $X^{(0)}:= X$,  $F_{X/k}^{(0)} := \mr{id}_X$, and $F_{X/k} := F^{(1)}_{X/k}$.

Let $\V$ be a  vector bundle  (i.e., a locally free coherent sheaf) on $X$.
We shall say that $\V$ is  {\bf decomposable} if there exist vector bundles $\V_1$ and $\V_2$ of positive ranks such that $\V\cong \V_1 \oplus \V_2$.
A vector bundle  is called {\bf indecomposable} if it is not decomposable.

According to ~\cite[Theorem 2]{Ati3},
any such vector bundle $\V$ is isomorphic to a direct sum of indecomposable vector bundles, unique up to reordering.
We refer to a vector bundle $\V'$ as an {\bf indecomposable component} of $\V$ if it  appears as  the direct summand in such a decomposition  (cf. ~\cite[Definition 2.3]{BiSu}).

\LSP
\subsection{Connections on vector bundles} \label{SS39}

A {\bf ($k$-)connection} on an $\mcO_X$-module $\V$ is defined as a $k$-linear morphism $\nabla : \V \rightarrow \Omega_X \otimes \V$ satisfying the Leibniz rule, i.e, $\nabla (a\cdot v) = a \cdot \nabla (v) + da  \otimes v$ for any local sections $a \in \mcO_X$, $v \in \V$.
A {\bf flat vector bundle} $X$ is a pair $(\V, \nabla)$  consisting of a vector bundle $\V$ on $X$ and a $k$-connection $\nabla$ on $\V$.
When $\V$ has rank one, such a pair is called a {\bf flat line bundle}.

If we are given two flat vector bundles $(\V, \nabla)$, $(\V', \nabla')$ on $X$, then the tensor product $\V \otimes \V' \left(:= \V \otimes_{\mcO_X} \V \right)$ admits a connection
\begin{align} \label{Eq285}
\nabla \otimes \nabla' : \V \otimes \V' \rightarrow \Omega_X \otimes (\V \otimes \V')
\end{align} 
 defined by  $(\nabla \otimes \nabla') (a \otimes b) = \nabla (a) \otimes b + a \otimes \nabla' (b)$ for any local sections $a \in \V$, $a' \in \V'$.
Also, given a flat vector bundle $(\V, \nabla)$, one can define a connection $\nabla^\vee$ on the dual bundle $\V^\vee := \mcH om_{\mcO_X} (\V, \mcO_X)$ in a natural manner (cf.~\cite[\S\,1.1]{Kat}).

\LSP
\subsection{$\mcD^{(\M)}$-modules} \label{SS32}

For each $\M \in \mbZ_{\geq 0}$,
let $\mcD_{X}^{(\M)}$ denote 
the sheaf of differential operators of level $\M$ on $X/k$ (cf. ~\cite[\S\,2.2]{PBer1}), where $\mr{Spec}(k)$ is equipped with the trivial $\M$-PD structure.
It defines a sheaf of (noncommutative) $k$-algebras and acts naturally on $\mcO_X$.
When there is no fear of confusion, we will write $\mcD^{(\M)} := \mcD_{X}^{(\M)}$ for simplicity.
There exists a natural inductive system of sheaves
  \begin{align} \label{Eqq2}
  \mcD^{(0)} \rightarrow \mcD^{(1)} \rightarrow \cdots \rightarrow \mcD^{(\M)} \rightarrow \cdots
  \end{align}
  (cf. ~\cite[\S\,2,2]{PBer1}).
The limit of this system 
   $\mcD^{(\infty)} :=  \varinjlim_{m \in \mbZ_{\geq 0}} \mcD^{(\M)}$ is called the sheaf of differential operators of {\it level $\infty$} and coincides with the sheaf discussed in ~\cite[\S\,16.8]{EGA4}.
   (Note that if $k$ has characteristic zero, then $\mcD_{X}^{(1)}=\mcD_{X}^{(2)}=\dots=\mcD_{X}$.)
 
For $\M, j  \in \mbZ_{\geq 0} \sqcup \{ \infty \}$, we write $\mcD_{\leq j}^{(\M)}$
 for the subsheaf of $\mcD^{(\M)}$ consisting of differential operators of order $\leq j$, with the convention that  $\mcD^{(\M)}_{\leq \infty} := \mcD^{(\M)}$.
 Thus, we have $\mcD^{(\M)} = \bigcup_{j \in \mbZ_{\geq 0}} \mcD^{(\M)}_{\leq j}$.
 We denote by  ${^L \mcD}^{(\M)}_{\leq j}$ (resp., ${^R \mcD}^{(\M)}_{\leq j}$)  the sheaf $\mcD^{(\M)}_{\leq j}$ endowed with a structure of  $\mcO_X$-module  induced by   left (resp., right) multiplication. 
Given an $\mcO_X$-module $\V$, we equip the tensor product $\mcD_{\leq}^{(\M)} \otimes \V := {^R}\mcD_{\leq j}^{(\M)} \otimes \V$ with the $\mcO_X$-module structure given by left multiplication.

Let $\mcE nd_k(\V)$ denote the sheaf of locally defined $k$-linear endomorphisms of $\mcE$ endowed with a structure of $\mcO_X$-module given by left multiplication. Given an $\mcO_{X}$-module $\V$, we equip the tensor product $\mcD_{\leq j}^{(m)} \otimes \V:={}^R \mcD_{\le j}^{(\M)} \otimes \V$ with the $\mcO_{X}$-module structure given by left multiplication.
A {\bf $\mcD^{(\M)}$-module structure} on $\mcE$ is an $\mcO_X$-linear morphism of sheaves of $k$-algebras $\nabla :  {^L}\mcD^{(\M)} \rightarrow \mcE nd_k (\mcE)$. Equivalently, specifying a $\mcD^{(\M)}$-module structure on $\mcE$ is a homomorphism  $\nabla$ extending the $\mcO_X$-module structure to a left $\mcD^{(m)}$-module structure.
By a {\bf $\mcD^{(\M)}$-bundle}, we mean a pair $(\V, \nabla)$ consisting of a vector bundle $\V$ on $X$ and a $\mcD^{(\M)}$-module structure $\nabla$ on it. A {\bf morphism of $\mcD^{(\M)}$-bundles} from $(\V,\nabla)$ to $(\V^{\prime},\nabla^{\prime})$ is defined as an $\mcO_{X}$-linear morphism $f:\V \rightarrow \V^{\prime}$ that is compatible with the respective $\mcD^{(\M)}$-module structures.

A {\bf  line  $\mcD^{(\M)}$-bundle} is  a $\mcD^{(\M)}$-bundle $(\V, \nabla)$ such that $\V$ has rank $1$.
For example, if 
 \begin{align} \label{eqttu}
 \nabla_{\mr{triv}}^{(\M)} :  {^L}\mcD^{(\M)} \rightarrow \mcE nd_k (\mcO_X)
 \end{align}
  denotes  the natural  $\mcD^{(\M)}$-action on $\mcO_X$  mentioned   above,
 then $(\mcO_X, \nabla_{\mr{triv}}^{(\M)})$ specifies 
 a line  $\mcD^{(\M)}$-bundle.

Let $(\V, \nabla)$ be a  $\mcD^{(\M)}$-bundle.
A {\bf  (line) $\mcD^{(\M)}$-subbundle} of $(\V, \nabla)$ is  
a (line) subbundle of $\V$ that is  stable under the $\mcD^{(\M)}$-action $\nabla$.
Each  $\mcD^{(\M)}$-subbundle $\V'$ of $\V$ naturally inherits 
a $\mcD^{(\M)}$-bundle structure, given by restricting $\nabla$.

Finally, it is verified (cf. ~\cite[Theorem 4.8]{BeOg}) that giving a $\mcD^{(0)}$-module structure on a fixed $\mcO_X$-module $\V$ amounts to giving a connection on $\V$.

\SSP
\begin{rem} \label{Reeem441}
Denote by $\mr{Str}  (X)$ the  abelian $k$-linear rigid tensor  category consisting  of $\mcD^{(\infty)}$-bundles.
According to ~\cite[Theorem 1.3]{Gie},
$\mr{Str}  (X)$ is equivalent to the category of $F$-divided sheaves on $X$.
Here, recall (cf. ~\cite{dSa}) that an {\bf $F$-divided sheaf} on $X$ is a collection
 \begin{align} \label{eq26}
 \{ (\V_\M, \alpha_\M) \}_{\M \in \mbZ_{\geq 0}},
 \end{align}
  where each pair $(\V_\M, \alpha_\M)$ consists of a vector bundle $\V_\M$ on $X^{(\M)}$ and  an $\mcO_{X^{(\M)}}$-linear isomorphism $\alpha_\M : F_{X^{(\M)}/k}^{(1)*} (\V_{\M +1}) \xrightarrow{\sim} \V_\M$.
 
Once we choose a $k$-rational point $x : \mr{Spec}(k) \rightarrow X$,
the assignment $(\V, \nabla) \mapsto x^* (\V)$ defines a tensor functor $\omega_x : \mr{Str} (X) \rightarrow \mr{Vec}_k$, where $\mr{Vec}_k$ denotes the category of finite-dimensional $k$-vector spaces.
The category $\mr{Str} (X)$ together with the functor $\omega_x$ forms a neutral Tannakian category (cf. ~\cite[\S\,2.2]{dSa}).
In particular, there exists a pro-algebraic group $k$-scheme $\pi_1^{\mr{str}} (X, x)$ (or simply $\pi_1^{\mr{str}} (X)$)  such that $\omega_x$ induces an equivalence
\begin{align} \label{Eq2229}
\mr{Str} (X) \cong \mr{Rep}_k (\pi_1^{\mr{str}}(X, x))
\end{align}
 between $\mr{Str} (X)$ and the category $\mr{Rep}_k (\pi_1^{\mr{str}}(X, x))$ of finite-dimensional $k$-representations of $\pi_1^{\mr{str}} (X, x)$.
The pro-algebraic group $\pi_1^{\mr{str}} (X, x)$  will be referred to as the {\bf stratified fundamental group} of $X$ (with respect to  the base point $x$).
\end{rem}

\LSP
\subsection{Frobenius pull-back of  $\mcD^{(\M)}$-modules} \label{SS32}

Suppose that $\M \in \mbZ_{\geq 0}$, and let
  $(\V, \nabla)$ be
 a $\mcD^{(\M)}$-bundle. 
  Note that the  natural 
short exact sequence  $0 \rightarrow \mcD^{(\M)}_{\leq p^{\M +1}-1} \rightarrow \mcD^{(\M)}_{\leq p^{\M +1}} \rightarrow \mcT^{\otimes p^{\M +1}}_X\rightarrow 0$ admits a canonical split injection $\mcT^{\otimes p^{\M +1}}_X \hookrightarrow \mcD^{(\M)}_{\leq p^{\M +1}}$  (cf.  ~\cite[\S\,2.5]{Wak12}).
The composite
  \begin{align} \label{Eq1030}
  F_{X/k}^{(\M +1) *}(\mcT_{X^{(\M +1)}})  \left(= \mcT^{\otimes p^{\M +1}}_X  \right) \rightarrow {^L}\mcD^{(\M)} \xrightarrow{\nabla} \mcE nd_{k} (\V),
\end{align}
of this injection with $\nabla$ factors through the inclusion $\mcE nd_{\mcO_X} (\V) \hookrightarrow \mcE nd_k (\V)$.
The resulting morphism
\begin{align} \label{eq2}
\psi^{\M +1} (\nabla) :  F_{X/k}^{(\M +1) *}(\mcT_{X^{(\M +1)}}) \rightarrow \mcE nd_{\mcO_X} (\V)
\end{align}
 is called the {\bf $p^{\M +1}$-curvature} of $\nabla$ (or, of $(\V, \nabla)$).

 This description of $p^{\M +1}$-curvature agrees with  the definition discussed  ~\cite{GLQ};
 see also  ~\cite[Definition 3.1.1]{LSQ} for a slightly different formulation  of a higher-level generalization of $p$-curvature.
When $\M =0$,  this definition  coincides with the classical definition of $p$-curvature (cf. ~\cite[\S\,5]{Kat}).
That is, the $p$-curvature  $\psi^1 (\nabla) : \left( F_{X/k}^* (\mcT_{X^{(1)}})=\right) \mcT_X^{\otimes p} \rightarrow \mcE nd_{\mcO_X} (\mcF)$ of a connection $\nabla$ is determined  by  $\partial^{\otimes p} \mapsto (\nabla_\partial)^p -\nabla_{\partial^p}$ for any local section $\partial \in \mcT_X$, where $\partial^p$ denotes the section of $\mcT_X$ corresponding to the $p$-th iterate $\partial \circ \cdots \circ \partial$ of $\partial$, and
write $\nabla_{(-)} := ((-) \otimes  \mr{id}_\V) \circ \nabla$.
 
Denote by  $\mcS ol (\nabla)$  the subsheaf of $\V$ on which $\mcD_+^{(\M)}$ acts by zero via $\nabla$, where $\mcD_+^{(\M)}$ denotes the kernel of the canonical projection $\mcD^{(\M)} \twoheadrightarrow \mcO_X$.
 This sheaf can be regarded as an $\mcO_{X^{(\M +1)}}$-submodule of $F_{X/k*}^{(\M +1)} (\V)$ via  the underlying homeomorphism of $F_{X/k}^{(\M +1)}$.

For a vector bundle $\mcG$ on $X^{(\M +1)}$ (or more generally, an $\mcO_{X^{(\M +1)}}$-module),
 there exists a canonical
 $\mcD^{(\M)}$-module structure 
\begin{align} \label{E445}
\nabla_{\mcG, \mr{can}}^{(\M)} : {^L}\mcD^{(\M)} \rightarrow\mcE nd_k (F_{X/k}^{(\M +1)*}(\mcG))
\end{align}
on the pull-back $F_{X/k}^{(\M +1)*}(\mcG)$  of $\mcG$
along  $F_{X/k}^{(\M +1)}$.
The resulting assignments  $\mcG \mapsto (F_{X/k}^{(\M +1) *}(\mcG), \nabla_{\mcG, \mr{can}}^{(\M)})$
and $(\V, \nabla) \mapsto \mcS ol (\nabla)$ together 
 define an equivalence of categories
\begin{align} \label{Eq577}
\left(\begin{matrix} \text{the category of} \\ \text{vector bundles on $X^{(\M +1)}$} \end{matrix} \right) \xrightarrow{\sim}
\left(\begin{matrix} \text{the category of $\mcD^{(\M)}$-bundles} \\ \text{with vanishing $p^{\M +1}$-curvature} \end{matrix} \right),
\end{align}
which maps  each  subbundle of a given vector bundle on $X^{(\M +1)}$ to a $\mcD^{(\M)}$-subbundle of the corresponding $\mcD^{(\M)}$-bundle (cf. ~\cite[Corollary 3.2.4]{LSQ}).

Moreover, if $s$ is another nonnegative integer, then the pull-back by $F_{X/k}^{(s)}$ induces 
an equivalence of categories
\begin{align} \label{Eq50}
\left(\begin{matrix} \text{the category of} \\ \text{$\mcD_{X^{(s)}}^{(\M)}$-bundles} \end{matrix} \right) \xrightarrow{\sim}
\left(\begin{matrix} \text{the category of} \\ \text{$\mcD^{(\M + s)}_{X}$-bundles} \end{matrix} \right)
\end{align}
(cf. ~\cite[Th\'{e}or\`{e}me 2.3.6]{PBer2}), which  
preserves the formation of $\mcD$-subbundles in the above sense.
If $(\V, \nabla)$ is a $\mcD^{(\M)}_{X^{(s)}}$-bundle and $(F_{X/k}^{(s)*}(\V), F_{X/k}^{(s)*}(\nabla))$ denotes its image under the above  equivalence, then 
it follows from  ~\cite[Proposition 2.2.4]{PBer2} and  (the proof of) ~\cite[Proposition 3.6]{GLQ} that
the pull-back  of $\psi^{\M +1} (\nabla)$  via $F_{X/k}^{(s)}$ coincides with $\psi^{\M + s +1} (F_{X/k}^{(s)*}(\nabla))$ via the natural inclusion $F_{X/k}^{(s)*}(\mcE nd_{\mcO_{X^{(s)}}}(\V)) \subseteq\mcE nd_{\mcO_X} (F_{X/k}^{(s)*}(\V))$, i.e.,
\begin{align} \label{eq4}
F_{X/k}^{(s)*} (\psi^{\M +1} (\nabla)) = \psi^{\M + s +1} (F_{X/k}^{(s)*}(\nabla)).
\end{align}
In particular, \eqref{Eq50} restricts to an equivalence between 
their subcategories 
\begin{align} \label{erw23}
\left(\begin{matrix} \text{the category of $\mcD_{X^{(s)}}^{(\M)}$-bundles} \\ \text{with nilpotent (resp., vanishing)} 
\\
\text{$p^{\M +1}$-curvature} 
\end{matrix} \right) \xrightarrow{\sim}
\left(\begin{matrix} \text{the category of $\mcD^{(\M + s)}_{X}$-bundles} \\ \text{with nilpotent (resp., vanishing)} \\
 \text{$p^{\M + s +1}$-curvature}
 \end{matrix} \right).
\end{align}

\LSP
\subsection{Complete flags on $\mcD^{(\M)}$-bundles} \label{SS12}

Let $\M$ be an element of  $\mbZ_{\geq 0}\sqcup \{\infty \}$ and 
   $\V$ a vector bundle on $X$ of rank $r \in \mbZ_{> 0}$.
Recall that a {\it complete flag} on $\V$ is defined as an $r$-step  increasing filtration
  \begin{align}
  0 = \V^0 \subseteq \V^1 \subseteq \cdots \subseteq \V^r = \V
  \end{align}
   such that
 all the subquotients $\V^{j}/\V^{j-1}$  are line bundles.
  The following definition introduces   the $\mcD^{(\M)}$-module version of this notion.

\SSP
\bde \label{Def3}
Suppose that $\V$ is equipped with  a $\mcD^{(\M)}$-module structure $\nabla$.
A {\bf  complete flag} on $(\V, \nabla)$ is  
a complete flag  $\{ \V^j \}_{j=0}^r$ on  $\V$ such that each subbundle $\V^j$ ($j=0, \cdots, r$) is preserved by  $\nabla$.
\ede
\SSP

The equivalence of categories \eqref{Eq50} leads to  the  following result, which reduces the  existence problem of complete flags for $\M < \infty$ to the case  $\M =0$.

\SSP
\bpr \label{Prop3200}
Let $\M$, $s$ be nonnegative integers and $(\mcG, \nabla)$  a $\mcD^{(\M)}_{X^{(s)}}$-bundle,
 and denote by  $(F^{(s)*}_{X/k} (\mcG), F^{(s)*}_{X/k} (\nabla))$ the $\mcD_X^{(\M +s)}$-bundle
corresponding to $(\mcG, \nabla)$ via \eqref{Eq50}.
Then, $(\mcG, \nabla)$ admits a complete flag  if and only if $(F^{(s)*}_{X/k} (\mcG), F^{(s)*}_{X/k} (\nabla))$ admits a complete flag.  
\epr
\SSP

Furthermore, the following assertion is well-known at least for $\M =0$ and follows  immediately  from the case $\M =0$  by applying 
Proposition \ref{Prop3200}.

\SSP
\bpr \label{Prop32}
Let $\M$ be a nonnegative integer and $(\V, \nabla)$ a $\mcD^{(\M)}$-bundle such that $\psi^{\M +1}(\nabla)$ is nilpotent.
Then, $(\V, \nabla)$ admits a complete flag.
\epr
\begin{proof}
By the equivalence of categories \eqref{erw23},
the problem is reduced to the case $\M = 0$.
We  prove the assertion for $\M = 0$ by induction on  the rank of $\V$.
There is nothing to prove when $\V$ has rank $1$.
For the induction step, suppose   that   $\mr{rank}(\V) > 1$. 
It follows from ~\cite[Corollary 5.5]{Kat} that $(\V, \nabla)$ admits a nonzero flat subbundle $(\mcF, \nabla_{\mcF})$ with vanishing $p$-curvature.
By the equivalence of categories \eqref{Eq577} for $\M = 0$,
$(\mcF, \nabla_{\mcF})$ descends to  a vector bundle $\mcG$ on $X^{(1)}$ along $F_{X/k}$.
Since the underlying space $X$ is a smooth projective curve,  we can find  a line subbundle  of $\mcG$.
This  corresponds  to a flat line subbundle  $(\mcH, \nabla_{\mcH})$ of  $(\V, \nabla)$ via  \eqref{Eq577}.
Let  $\overline{\nabla}$ denote the connection  on  $\overline{\V} := \V/\mcH$ induced from $\nabla$ via  the quotient $\pi : \V \twoheadrightarrow \overline{\V}$. Since  $\overline{\nabla}$ has nilpotent $p$-curvature,  the induction hypothesis  implies that  $(\overline{\V}, \overline{\nabla})$ admits  a  complete flag   $\{ \overline{\V}^j \}_j$.
Hence, the collection  $\{ \V^j \}_{j}$ defined as  $\V^0 := 0$ and $\V^{j+1} := \pi^{-1} (\V^{j})$ ($j \in \mbZ_{\geq 0}$) forms the  desired complete  flag on $(\V, \nabla)$.
This completes the proof.
\end{proof}

\LSP
\section{The case of projective line} \label{S234}
\LSP

This section focuses on  the case where the underlying curve $X$ has genus $0$, i.e., $X$ is isomorphic to the projective line.
In particular, we describe  vector bundles on $X$ admitting a $\mcD^{(\M)}$-module structure  (cf. Proposition \ref{Prop6}) and prove the existence of complete flags (cf. Theorem \ref{Th4}).

\LSP
\subsection{Vector bundles admitting a $\mcD^{(\M)}$-module structure} \label{SS46}

Denote by $\mbP$ the projective line over $k$ throughout  this section.
(We use the notation ``$\mbP$" as opposed to the usual notation ``$\mbP^1$" because later on we will need to consider its $\M$-th Frobenius twist for $\M \in \mbZ_{\geq 0}$ and $\mbP^{(\M)}$ is notationally and typographically simpler than $(\mbP^1)^{(\M)}$.)
 For each $\ell \in \mbZ$, the unique (up to isomorphism) line bundle on $\mbP$ of degree $\ell$ is denoted by $\mcO_\mbP$.
Recall the Birkhoff-Grothendieck theorem (cf.  ~\cite{Grot}), which states that any vector bundle  on $\mbP$ is isomorphic to a direct sum of  line bundles of the form $\mcO_{\mbP}$.
In particular,   
 every  indecomposable vector bundle on $\mbP$ has rank $1$.
 Applying this theorem, one can  prove the following result, generalizing 
  the classification of flat vector bundles  in the genus-$0$ case (cf.  ~\cite{Ati1}, ~\cite{Bis1}, ~\cite{BiSu}, and ~\cite{Wei}).

\SSP
\bpr \label{Prop6}
Let $\M$ be a nonnegative integer  and $\V$  a vector bundle on $X \left(= \mbP \right)$ of rank $r \in \mbZ_{> 0}$.
Then, the following three conditions (a)-(c) are equivalent:
\begin{itemize}
\item[(a)]
$\V$ admits a $\mcD^{(\M)}$-module structure;
\item[(b)]
$\V$ admits a $\mcD^{(\M)}$-module structure with vanishing $p^{\M+1}$-curvature;
\item[(c)]
$\V \cong \bigoplus\limits_{i=1}^r \mcO_{\mbP} (p^{\M +1}_i)$ for some $\ell_1, \cdots, \ell_r \in \mbZ$.
\end{itemize}
\epr
\begin{proof}
(b) $\Rightarrow$ (a) is clear.
To prove  (a) $\Rightarrow$ (c),
we assume  that $\V$ admits a $\mcD^{(\M)}$-module structure $\nabla_\V$.
By the equivalence of categories \eqref{Eq50}, $(\V, \nabla_\V)$ corresponds to a  flat vector bundle  $(\mcG, \nabla_\mcG)$ on $X^{(\M)} \left(= \mbP^{(\M)} \right)$ via pull-back by $F^{(\M)}_{\mbP/k}$.
It follows from the Birkhoff-Grothendieck theorem that  $\mcG \cong \bigoplus_{i=1}^r \mcO_{\mbP^{(\M)}} (\ell'_i)$ for some $\ell'_1, \cdots, \ell'_r \in \mbZ$.
As proved in the references listed above,
there exists a collection $(\ell_1, \cdots, \ell_r) \in \mbZ^r$ with
$\ell'_i = p \ell_i$   for every $i$.
Hence, we have 
\begin{equation*}
\V \cong F^{(\M)*}_{{\mbP}/k} (\mcG) \cong F^{(\M)*}_{{\mbP}/k} (\bigoplus_{i=1}^r \mcO_{\mbP^{(\M)}} (p \ell_i)) \cong \bigoplus_{i=1}^r F^{(\M) *}_{{\mbP}/k} (\mcO_{\mbP^{(\M)}} (p\ell_i)) \cong \bigoplus_{i=1}^r \mcO_{\mbP} (p^{\M +1} \ell_i).
\end{equation*}
Thus, $\V$ satisfies 
the condition (c).

Finally, we shall prove (c) $\Rightarrow$ (b).
Assume that $\V  =  \bigoplus\limits_{i=1}^r \mcO_{\mbP} (p^{\M +1}\ell_i)$ for some  $\ell_1, \cdots, \ell_r \in \mbZ$.
By putting  $\mcG := \bigoplus_{i=1}^r \mcO_{\mbP^{(\M +1)}} (\ell_i)$,  we obtain a chain  of isomorphisms
\begin{align}
F^{(\M +1) *}_{{\mbP}/k} (\mcG) \xrightarrow{\sim}\bigoplus_{i=1}^r F^{(\M +1)*}_{{\mbP}/k} (\mcO_{\mbP^{(\M +1)}} (\ell_i))  \xrightarrow{\sim}\V.
\end{align}
Hence,
$\nabla_{\mcG, \mr{can}}^{(\M)}$ (cf. \eqref{E445}) defines a  $\mcD^{(\M)}$-module structure on $\V$ (with vanishing $p^{\M +1}$-curvature) via this composite.
In particular, $\V$ satisfies the condition (b).
This completes the proof.
\end{proof}

\LSP
\subsection{The existence of complete flags} \label{SS49}

Applying Proposition \ref{Prop6}, we can prove the following assertion,
which is  a part of Theorem \ref{ThA}.

\SSP
\bt \label{Th4}
Suppose that $X=\mbP$.
Let $\M$ be an element of $\mbZ_{\geq 0} \sqcup \{ \infty \}$ and $(\V, \nabla)$  a (nonzero) $\mcD^{(\M)}$-bundle on $X$.
Then, $(\V, \nabla)$ admits  a complete flag.
\et
\begin{proof}
In the case of  $\M = \infty$, we recall from ~\cite[Theorem 2.2]{Gie}  that the stratified fundamental group $\pi_1^\mr{str} (\mbP)$  of $\mbP$  (cf. Remark \ref{Reeem441})   is trivial.
That is, every $\mcD^{(\infty)}$-bundle is isomorphic to a direct sum of finitely many copies of $(\mcO_{\mbP}, \nabla_\mr{triv}^{(\infty)})$, so the assertion for $\M = \infty$ follows immediately from this fact.

For $\M \in \mbZ_{\geq 0}$,
the problem is reduced to the case of $\M = 0$ by Proposition \ref{Prop3200}.
Let us prove the assertion for $\M  = 0$
by  induction  on $r := \mr{rank} (\V)$.
There is nothing to prove in the base case of rank $r =1$.
To consider  the induction step, we assume that the assertion
for rank $r-1$ flat vector bundles
 (where $r \geq 2$) has been proved.
The equivalence (a) $\Leftrightarrow$ (c) in Proposition \ref{Prop6} implies
\begin{align} \label{Eq57}
\V \cong  \bigoplus_{i=1}^s  \mcO_{\mbP} (p\ell_i)^{\oplus e_i}
\end{align}
 for  $\ell_1 >  \ell_2 > \cdots  > \ell_s$ ($s \in \mbZ_{>0}$) and $e_1, \cdots, e_s \in \mbZ_{>0}$  (with $\sum_{i=1}^s e_i = r$).
We write
 $\mcG:= \bigoplus_{i=1}^s  \mcO_{\mbP} (\ell_i)^{\oplus e_i}$.
Then, there exists a natural isomorphism $\alpha : \V \xrightarrow{\sim} F^{(1) *}_{\mbP/k}(\mcG)$, which allows us to regard $\nabla_{\mcG, \mr{can}}^{(0)}$ as  a connection on   $\V$.
Since the difference $\nabla - \nabla_{\mcG, \mr{can}}^{(0)}$ is $\mcO_{\mbP}$-linear, 
 $\nabla$ can be expressed as 
\begin{align} \label{Eq61}
\nabla= \nabla_{\mcG, \mr{can}}^{(0)} + \sum_{1 \leq i, j \leq s}h_{i, j},
\end{align}
where $h_{i, j} \in \mr{Hom}_{\mcO_{\mbP}} (\mcO_{\mbP} (p\ell_i)^{\oplus e_i}, \Omega_{\mbP} \otimes (\mcO_{\mbP} (p\ell_j)^{\oplus e_j}))$.
Note that  $h_{i, j} = 0$ whenever
$i < j$, in particular, that
 $h_{1, j} = 0$ for every $j=2, \cdots, s$.
Since $\mr{End}_{\mcO_{\mbP}} (\mcO_{\mbP} (p\ell_1)) = k$,
the endomorphism $h_{1, 1}$ can be given by an $e_1 \times e_1$ matrix $A \in M_{e_1} (k)$.
Choose an automorphism $\eta$ of $\mcO_{\mbP} (p \ell_1)^{\oplus e_1}$ such that, if $P \in M_{e_1} (k)$  is the corresponding  matrix, then  $P^{-1}A P$ is upper triangular.
The gauge transformation by the automorphism $\widetilde{\eta} := \eta \oplus \bigoplus_{i=2}^s \mr{id}_{\mcO_{\mbP} (p\ell_i)^{\oplus e_i}}$ transforms  $\nabla$ into another  connection 
$\nabla'$.
The $1$-st direct summand $\mcL'$ in  $\mcO_{\mbP} (p\ell_1)^{\oplus e_1}$, regarded as a line subbundle of $\V$ via $\alpha$, satisfies   $\nabla' (\mcL') \subseteq \Omega_{\mbP} \otimes \mcL'$.
 It follows that the line subbundle $\mcL$ of $\V$ corresponding to $\mcL'$ via $\widetilde{\eta}$ satisfies  $\nabla (\mcL) \subseteq \Omega_{\mbP^{1}} \otimes \mcL$.
Write $\overline{\V} := \V/\mcL$, and let 
$\overline{\nabla}$ be the connection on $\overline{\V}$ induced from $\nabla$ via  the natural quotient $\pi : \V \twoheadrightarrow \overline{\V}$.
By the induction hypothesis,  $(\overline{\V}, \overline{\nabla})$ admits a complete flag $\{ \overline{\V}^j \}_{j=0}^{r-1}$.
For each $j =0, \cdots, r$, we shall set  $\V^0 := 0$ and $\V^j := \pi^{-1} (\overline{\V}^{j-1})$ ($j = 1, \cdots, r$).
Then, the resulting filtration $\{ \V^j \}_j$ forms the desired complete flag.
This completes the proof.
\end{proof}

\LSP
\section{The case of elliptic curves} \label{S45}
\LSP

This section  is devoted to studying the case of curves of genus $1$.
Our discussion relies on  Atiyah's classification of indecomposable vector bundles on such curves.
As a consequence, we establish a part of our main theorem, proving  the existence of complete flags  for genus-$1$ curves  (cf. Theorem \ref{Th2}).

Throughout  this section, this $g=1$ assumption will be in force.
In particular, by our assumption, $\Omega_X$ is trivialized.
We fix an identification $\mcO_{X}=\Omega_{X}$ by choosing a global generator of $\Omega_{X}$.
Additionally, we fix a $k$-rational point $\sigma$ of $X$, regarded as an effective divisor of degree one.
 
\LSP
\subsection{Atiyah's classification of indecomposable bundles} \label{SS12}

We begin with reviewing  the study of indecomposable vector bundles on a genus-$1$ curve, discussed in ~\cite{Ati1}.

Let $V$ be a finite-dimensional $k$-vector space and $\mcE$ a vector bundle on $X$.
There exists a sequence of natural isomorphisms
\begin{align} \label{Eq223}
\mr{Ext}^1_{\mcO_X} (\mcE, V \otimes_k \mcL) \xrightarrow{\sim} \mr{Hom}_{k} (\mr{Ext}^1_{\mcO_X}(\mcE, \mcL)^\vee, V) \xrightarrow{\sim} \mr{Hom}_k (H^0 (X, \mcL^\vee \otimes \mcE), V),
\end{align}
where the second isomorphism  arises from  Serre duality and our fixed identification  $\mcO_X = \Omega_X$.
Hence, each extension $\varepsilon : 0 \rightarrow V \otimes_k \mcL \rightarrow \mcF \rightarrow \mcE \rightarrow 0$ corresponds to a $k$-linear morphism
\begin{align} \label{Eq444}
\mr{ext}_\mcF \ \left(\text{or} \ \mr{ext}_\varepsilon \right) : H^0 (X, \mcL^\vee \otimes \mcE) \rightarrow V
\end{align}
via  \eqref{Eq223}.
If $V = H^0 (X, \mcL^\vee \otimes \mcE)$,
then $\mr{ext}_\mcF$ defines a $k$-linear endomorphism of $H^0 (X, \mcL^\vee \otimes \mcE)$.

Now, recall  from  ~\cite[Theorem 5, (i)]{Ati1} that, for each positive integer $r$,
there exists a unique (up to isomorphism) indecomposable vector bundle $\mcF_r$ on $X$ of rank $r$ and degree $0$ with $H^0 (X, \mcF_r) \neq 0$.
Moreover, $\mr{dim}_k (H^0 (X, \mcF_r)) =1$  and the cokernel of the unique (up to multiplication by an element of  $k^\times$) injection $\mcO_X \hookrightarrow \mcF_r$  is isomorphic to $\mcF_{r -1}$ (with the convention that  $\mcF_0 := 0$).
 In other words,
 by fixing  an isomorphism $\mr{Coker}(\mcO_X \hookrightarrow \mcF_r) \xrightarrow{\sim} \mcF_{r-1}$, 
  we obtain  a short exact sequence
 \begin{align} \label{Eq700}
 0 \rightarrow \mcO_X \rightarrow \mcF_{r} \rightarrow \mcF_{r-1} \rightarrow 0
 \end{align}
 (cf. ~\cite[Lemma 15, (i), and Theorem 5, (i)]{Ati1}).
This extension corresponds to the identity morphism of $H^0 (X, \mcF_{r-1})$ via \eqref{Eq223} in the case where $V = H^0 (X, \mcF_{r-1}) \left(\cong k \right)$ and $\mcL = \mcO_X$.

Denote by
\begin{align} \label{eqwhe}
\mr{Fil}^\bullet (\mcF_r) := \{ \mr{Fil}^j (\mcF_r) \}_{j=0}^r
\end{align}
the filtration on $\mcF_r$ defined in such a way that for each $j =0, \cdots, r$,  the subsheaf  $\mr{Fil}^j (\mcF_r)$ is  the kernel of the composite surjection
$\mcF_r \twoheadrightarrow  \mcF_{r-1}\twoheadrightarrow \cdots \twoheadrightarrow \mcF_{r-j}$.
Then, we have $\mr{Fil}^j (\mcF_r)/\mr{Fil}^{j-1} (\mcF_r) \cong \mcO_X$ for every $j =1, \cdots, r$.

Next,take  an  arbitrary indecomposable  vector bundle $\V$ on $X$ of rank $r \in \mbZ_{> 0}$ and degree $d \in \mbZ$.
Define $h := \mr{gcd} (r, d)$.
Following ~\cite[Theorem 6]{Ati1}, we shall construct, by induction on $j$, proper subbundles $\mr{Fil}^j (\V)$ ($j=0, 1, 2, \cdots$) of $\mcE$, starting with $\mr{Fil}^0 (\V) := 0$ and satisfying the following conditions:
\begin{itemize}
\item
 $\mcE/\mr{Fil}^j (\mcE)$ is  indecomposable;
\item
The equality $h = \mr{gcd}(r_j, d_j)$ holds, where $r_j := \mr{rank}(\V/\mr{Fil}^j (\V))$ and $d_j := \mr{deg}(\V/\mr{Fil}^j (\V))$ (hence $(r_0, d_0) = (r, d)$).
\end{itemize}

To construct $\mr{Fil}^j (\mcE)$, assume that a proper subbundle $\mr{Fil}^{j-1} (\V)$ ($j \geq 1$) has been constructed, and that $r_{j-1} \nmid d_{j-1}$.
Under these assumptions,
the vector bundle $\mcL_j^\vee \otimes (\V/\mr{Fil}^{j-1} (\V))$, where $\mcL_j := \mcO_X (\left\lfloor \frac{d_{j-1}}{r_{j-1}}\right\rfloor \sigma)$,  has positive degree strictly smaller than $r_{j-1}$.
It follows that the natural morphism
\begin{align}
\alpha_j : H^0 (X, \mcL^\vee_j \otimes (\V/\mr{Fil}^{j-1}(\V))) \otimes_k \mcL_j \rightarrow \V/\mr{Fil}^{j-1}(\V)
\end{align}
is nonzero and  injective,  and its cokernel $\mr{Coker}(\alpha_j)$ defines   an indecomposable vector bundle (cf. ~\cite[Lemma 15, (ii)]{Ati1}).
We define  $\mr{Fil}^j (\V)$ as  the inverse image of $\mr{Im}(\alpha_j)$ along the quotient $\V \twoheadrightarrow \V/\mr{Fil}^{j-1} (\V)$.
Since  $\V/\mr{Fil}^j (\V)$ coincides with $\mr{Coker} (\alpha_j)$, it is indecomposable, and we obtain  a short exact sequence
\begin{align} \label{Eq334}
0 \rightarrow H^0 (X, \mcL_j^\vee \otimes (\V/\mr{Fil}^{j-1}(\V))) \otimes_k \mcL_j \rightarrow \V/\mr{Fil}^{j-1} (\V) \rightarrow \V/\mr{Fil}^{j} (\V) \rightarrow 0.
\end{align}
Equivalently, we obtain an isomorphism of vector bundles
\begin{align} \label{Eq456}
H^0 (X, \mcL_j^\vee \otimes (\V/\mr{Fil}^{j-1}(\V))) \otimes_k \mcL_j \xrightarrow{\sim} \mr{Fil}^{j}(\V)/\mr{Fil}^{j-1}(\V) \left(\neq 0 \right).
\end{align}
Moreover, the equality
\begin{align} \label{Eq3333}
\mr{dim}_k (H^0 (X, \mcL_j^\vee \otimes (\V/\mr{Fil}^{j-1}(\V)))) = \mr{dim}_k (H^0 (X, \mcL_j^\vee \otimes (\V/\mr{Fil}^j (\V))))
\end{align}
holds (cf. ~\cite[Lemma 15, (ii)]{Ati1}).
Under a suitable identification $H^0 (X, \mcL_j^\vee \otimes (\V/\mr{Fil}^{j-1}(\V))) = H^0 (X, \mcL_j^\vee \otimes (\V/\mr{Fil}^j (\V)))$, the extension \eqref{Eq334} corresponds to the identity morphism of $H^0 (X, \mcL_j^\vee \otimes (\V/\mr{Fil}^j (\V)))$ via \eqref{Eq223}.
This means  $\mr{ext}_{\V/\mr{Fil}^{j-1}(\V)} = \mr{id}_{H^0 (X, \mcL_j^\vee \otimes (\V/\mr{Fil}^j (\V)))}$ in the notation of \eqref{Eq444}.
Since
\begin{align}
r_j = r_{j-1}-\left(d_{j-1}-r_{j-1}\cdot \left\lfloor \frac{d_{j-1}}{r_{j-1}}\right\rfloor \right),
\hspace{5mm}
 d_j = d_{j-1}-\left\lfloor \frac{d_{j-1}}{r_{j-1}}\right\rfloor \cdot \left(d_{j-1}-r_{j-1}\cdot \left\lfloor \frac{d_{j-1}}{r_{j-1}}\right\rfloor \right),
\end{align}
we have $\mr{g.c.d}(r_j, d_j) = \mr{g.c.d}(r_{j-1}, d_{j-1}) = h$ and $r_j < r_{j-1}$, $d_j < d_{j-1}$, ensuring that the process terminates.

By repeating this construction, we end up with an increasing filtration
$\{ \mr{Fil}^j (\V) \}_{j=0}^m$
of $\V$ for some  $m \in \mbZ_{> 0}$ 
such that  $r_m \mid d_m$ (hence $r_m = h$).
According to ~\cite[Theorem 5, (ii)]{Ati1},
there exists a line bundle $\mcL_{m+1}$ of degree $\frac{d_m}{r_m}$ along with  an isomorphism  $\alpha_\V : \V/\mr{Fil}^m (\mcE) \xrightarrow{\sim} \mcL_{m+1}\otimes \mcF_{r_m}$.
Using this, we extend  $\{ \mr{Fil}^j (\V) \}_{j=0}^m$ to a new  filtration
\begin{align} \label{Eq3334}
\mr{Fil}^\bullet (\V) := \{ \mr{Fil}^j (\V) \}_{j=0}^{\ell},
\end{align}
where $\ell := m + r_m$, by defining  
\begin{align}
\mr{Fil}^j(\V) := \mr{Ker} \left(\V \twoheadrightarrow \V/\mr{Fil}^{m}(\V) \xrightarrow{\alpha_\V} \mcL_{m+1} \otimes  \mcF_{r_m} \twoheadrightarrow \mcL_{m+1} \otimes  (\mcF_{r_m}/\mr{Fil}^{j-m} (\mcF_{r_m})) \right)
\end{align}
for $j=m +1, \cdots, \ell$ (hence $\mr{Fil}^\ell (\V) = \V$).
Setting  $\mcL_{j}:= \mcL_{m+1}$ for $j=m+1, \cdots, \ell$, we obtain a collection of line bundles
\begin{align} \label{Eq299}
L (\V) := (\mcL_1, \cdots, \mcL_{m (\mcE)}, \cdots,  \mcL_{\ell}),
\end{align}
where $m (\mcE) := m \left( \geq 0 \right)$.
Note that both \eqref{Eq456} and \eqref{Eq3333} hold even for $j \geq m+1$ (cf. ~\cite[Lemma 16]{Ati1}).
Additionally, if   $m (\mcE) > 0$, then we have  $\mcL_1 \cong \mcO_X (\left\lfloor \frac{r}{d}\right\rfloor \sigma))$.
Furthermore,  in the special case where  $\V= \mcF_r$ ($r \in \mbZ_{\geq 0}$),
 the filtration \eqref{Eq3334} coincides with the filtration \eqref{eqwhe}.

\LSP
\subsection{Connections on indecomposable bundles} \label{SS101}

As a preparatory  step for establishing   the  main result of this section, we first prove  several properties of connections on  indecomposable vector bundles.
Let us fix an indecomposable vector bundle $\V$ on $X$, and let  $\{ \mr{Fil}^j (\V) \}_{j=0}^\ell$  be the increasing  filtration on $\V$ constructed as in
\eqref{Eq3334}.
Also, write $L (\V) = (\mcL_1, \cdots, \mcL_{\ell})$.

\SSP
\bpr \label{Lem47}
Suppose that $\V$ admits a connection $\nabla$.
Then, for each $j=1, \cdots, \ell$, the subbundle $\mr{Fil}^j(\mcE)$ is stable under $\nabla$, i.e.,  $\nabla (\mr{Fil}^j (\V)) \subseteq \Omega_X \otimes \mr{Fil}^j (\V)$.
\epr
\begin{proof}
Suppose, on the contrary,  that the filtration $\{ \mr{Fil}^j (\V) \}_j$ is not preserved by  $\nabla$.
In particular, the integer 
\begin{align}
j_0 := \mr{min} \{ j \, | \, \nabla (\mr{Fil}^{j_0} (\mcE)) \not\subseteq \Omega_X \otimes \mr{Fil}^{j_0} (\mcE) \} \left(> 0 \right)
\end{align}
is well-defined.
We set  $\mcL := \mcL_{j_0}$, $\mcM := \V/\mr{Fil}^{j_0 -1}(\V)$,  
 and  $\mcN := \mr{Fil}^{j_0}(\V)/\mr{Fil}^{j_0 -1} (\V)$  ($\cong H^0 (X, \mcL^\vee_{j_0} \otimes \mcM) \otimes_k\mcL_{j_0}$ by \eqref{Eq456}).
Since $\mr{Fil}^{j_0 -1} (\V)$ is closed  under $\nabla$,
this connection induces a connection $\overline{\nabla}$ on $\mcM$.
It follows that there exists an injection
$\iota : \mcL\hookrightarrow \mcN \left(\subseteq \mcM \right)$ with 
 $\overline{\nabla} (\mr{Im}(\iota)) \not\subseteq \Omega_X \otimes \mcN$.
Let us consider the $\mcO_X$-linear  composite
\begin{align} \label{Eq288}
\mcD^{(0)}_{\leq 1} \otimes \mcL \hookrightarrow \mcD^{(0)} \otimes \mcM \xrightarrow{\overline{\nabla}} \mcM,
\end{align}
where the first arrow denotes  the tensor product of the natural inclusion $\mcD^{(0)}_{\leq 1} \hookrightarrow \mcD^{(0)}$ and $\iota$.
This composite  restricts to  $\iota : \mcL \left(= \mcD^{(0)}_{\leq 0} \otimes \mcL\right) \rightarrow \mcN$, and  induces 
 a {\it nonzero} $\mcO_X$-linear morphism 
 $\eta : \mcT_X \otimes \mcL \rightarrow \mcM/\mcN$ via the quotient 
 \begin{align}
 \mcD^{(0)}_{\leq 1} \otimes \mcL \twoheadrightarrow (\mcD_{\leq 1}^{(0)}/\mcD_{\leq 0}^{(0)}) \otimes \mcL = \mcT_X \otimes \mcL.
 \end{align}
Since  $\mcT_X \cong \mcO_X$,
the morphism
\begin{align} \label{Eq1111}
\eta' : \left(k \cong  \right)H^0 (X, \mcT_X) \rightarrow H^0 (X, \mcL^\vee \otimes (\mcM/\mcN))
\end{align}
induced by $\eta$ is nonzero.
Also, we obtain   the following commutative diagram:
\begin{align} \label{Eq125}
\vcenter{\xymatrix@C=46pt@R=36pt{
0 \ar[r] & \mcL\ar[d]^-{\iota}  \ar[r]^-{\mr{inclusion}} & \mcD^{(0)}_{\leq 1} \otimes \mcL \ar[r]^-{\mr{quotient}} \ar[d]^-{\eqref{Eq288}}  & \mcT_X \otimes \mcL \ar[r] \ar[d]_-{\wr}^-{\mr{id}} & 0\\
 0 \ar[r] & \mcN\ar[d]_-{\wr}^-{\mr{id}}  \ar[r]^-{\mr{inclusion}} & \pi^{-1} (\mcT_X \otimes \mcL) \ar[r]^-{\pi |_{ \pi^{-1} (\mcT_X \otimes \mcL)}} \ar[d]^-{\mr{inclusion}}  & \mcT_X \otimes \mcL \ar[r] \ar[d]^-{\eta} & 0\\
 0 \ar[r] & \mcN \ar[r]^-{\mr{inclusion}}  & \mcM \ar[r]^-{\pi} & \mcM/\mcN \ar[r] & 0, 
 }}
\end{align}
where $\pi$ denotes the natural quotient and  all  horizontal sequences are exact.
The extension defined by  the lower horizontal sequence in  \eqref{Eq125},  denoted by $\varepsilon$,  
corresponds to the identity morphism of $H^0 (X, \mcL^\vee \otimes (\mcM/\mcN))$ via \eqref{Eq223}, i.e., $\mr{ext}_\varepsilon = \mr{id}_{H^0 (X, \mcL^\vee \otimes (\mcM/\mcN))}$.
The extension class of the middle horizontal sequence is nontrivial, since it corresponds to $\eta'$ via  the composite isomorphism
 \begin{align}
 \mr{Ext}^1 (\mcT_X \otimes \mcL, \mcN)
 &\xrightarrow{\sim}\mr{Ext}^1 (\mcT_X \otimes \mcL, H^0 (X, \mcL^\vee \otimes \mcM)\otimes_k  \mcL) \\
  &\xrightarrow{\sim}\mr{Ext}^1 (\mcT_X \otimes \mcL, H^0 (X, \mcL^\vee \otimes (\mcM/\mcN))\otimes _k \mcL) \\
 &\xrightarrow{\sim}
 \mr{Hom}_k (H^0 (X, \mcT_X), H^0 (X, \mcL^\vee \otimes (\mcM/\mcN)))  \notag
 \end{align}
 (cf. \eqref{Eq223}, \eqref{Eq3333}).
On the other hand, 
the upper horizontal arrow in \eqref{Eq125} splits,  since  the fixed  global generator of $\mcT_X$ provides  a splitting of  the natural short exact sequence $0 \rightarrow \mcO_X \left(= \mcD^{(0)}_{\leq 0} \right) \rightarrow {^R}\mcD^{(0)}_{\leq 1} \rightarrow \mcT_X \rightarrow 0$.
The composite of the resulting splitting $\mcT_X \otimes \mcL \hookrightarrow \mcD^{(0)}_{\leq 1} \otimes \mcL$ and the inclusion $\mcD^{(0)}_{\leq 1} \otimes \mcL \hookrightarrow \pi^{-1} (\mcT_X \otimes \mcL)$ (i.e., the upper middle vertical arrow in \eqref{Eq125}) specifies a splitting of the middle horizontal sequence, leading to  a contradiction.
Thus,
 the original assumption that $\nabla$ does not preserve the filtration $\{ \mr{Fil}^j (\V)\}_j$ must be false,  completing  the proof.
\end{proof}
\SSP

\bpr \label{Lem45}
 $\V$ admits a connection  if and only if
$p \mid \mr{deg}(\mcL_j)$ for every $j=1, \cdots, \ell$.
\epr
\begin{proof}
First, we shall consider the ``only if" part of the required equivalence.
Suppose that there exists a connection $\nabla$ on $\V$. 
By Proposition \ref{Lem47},
$\nabla$ preserve the filtration $\{ \mr{Fil}^j (\V) \}_j$.
Hence, for each $j= 1, \cdots, \ell$,
 $\nabla$ induces a connection on
 $\mr{Gr}^j (\V) := \mr{Fil}^j (\V)/\mr{Fil}^{j-1} (\V)$.
 Since $\mr{Gr}^j (\V)$ is a direct sum of finitely many copies of $\mcL_j$,
the line bundle $\mcL_j$ admits a connection (cf. ~\cite[Lemma 2.4]{Bis1}).
 Then, it follows from ~\cite[Theorem 1.1]{Bis1} or ~\cite[Theorem 1.1]{BiSu} that $p \mid \mr{deg}(\mcL_j)$.
 
 Next, we shall consider the ``if" part.
 Assume that $p \mid \mr{deg}(\mcL_j)$ for every $j$.
 By this assumption, each $\mcL_j$ admits  a connection $\nabla_j$  (cf. ~\cite[Theorem 1.1]{Bis1} or ~\cite[Theorem 1.1]{BiSu}).
 We now construct  a connection on $\V/\mr{Fil}^j (\mcE)$ by descending  induction  on $j =0, \cdots, \ell -1$.
 In the case of $j=\ell-1$, we recall that $\V/\mr{Fil}^{\ell -1}(\V) \left(= \mr{Gr}^{\ell} (\V) \right)$ is isomorphic to a direct sum  of copies of $\mcL_\ell$.
 Hence, the direct sum  of $\nabla_\ell$ defines  a connection on $\V/\mr{Fil}^{\ell -1}(\V)$, establishing  the base step of our induction argument.
 To consider the induction step, we assume that
  a connection $\nabla'_j$ on $\mcE/\mr{Fil}^{j} (\mcE)$ has been constructed
  for a fixed integer $j \in \{1, \cdots, \ell -1 \}$.
 Recall that $\mr{Gr}^{j} (\mcE) \cong H^0 (X, \mcL_j^\vee \otimes (\V/\mr{Fil}^{j}(\V))) \otimes_k \mcL_j$ (cf. \eqref{Eq334} and \eqref{Eq3333}) and the natural extension
 \begin{align} \label{Eq203}
 \varepsilon : 0 \rightarrow \mr{Gr}^j (\V) \rightarrow \V/\mr{Fil}^{j-1} (\V) \rightarrow \V/\mr{Fil}^{j} (\V) \rightarrow 0
 \end{align}
 corresponds to the identity morphism of $H^0 (X, \mcL_j^\vee \otimes (\V/\mr{Fil}^j (\V)))$ via 
 \eqref{Eq223}.
 Consider the corresponding morphism 
 \begin{align}
 \alpha : H^0 (X, \mcL_j^\vee \otimes (\mcE/\mr{Fil}^{j} (\mcE))) \otimes_k \mcO_X \rightarrow \mcL_j^\vee \otimes (\mcE/\mr{Fil}^{j} (\mcE)).
\end{align}
 One can find a connection $\nabla''_{j}$ on $H^0 (X, \mcL_j^\vee \otimes (\mcE/\mr{Fil}^{j} (\mcE))) \otimes_k \mcO_X$ 
 such that, under the fixed identification $\mcO_X = \Omega_X$,  the endomorphisms  of $H^0 (X, \mcL_j^\vee \otimes (\mcE/\mr{Fil}^{j} (\mcE)))$  induced from $\nabla''_{j}$ and $\nabla_j^\vee \otimes \nabla'_j$ (cf. \eqref{Eq285}) coincide  via applying the functor $H^0 (X, -)$ to $\alpha$.
 The connection  $\nabla''_j \otimes \nabla_j$  can be transposed into a connection $\nabla'''_j$  on $\mr{Gr}^j (\V)$ by using  the isomorphism  $H^0 (X, \mcL_j^\vee \otimes (\V/\mr{Fil}^j (\V))) \otimes_k \mcL_j \xrightarrow{\sim} \mr{Gr}^j (\mcE)$ (cf.  \eqref{Eq456}).
 In particular, this yields  a connection $\nabla'''_j \otimes (\nabla'_j)^\vee$ on  $\mr{Gr}^j (\V) \otimes (\V/\mr{Fil}^j (\V))^\vee$.
 Denote by $\mcK^\bullet  [\nabla'''_j \otimes (\nabla'_j)^\vee]$ the complex of sheaves  defined to be $\nabla'''_j \otimes (\nabla'_j)^\vee$   concentrated at degrees $0$ and $1$.    
 The Hodge-to-de Rham spectral sequence associated to  this complex  gives
 an exact sequence
 \begin{align} \label{Eq4421}
 \mbH^1 (X, \mcK^\bullet [\nabla'''_j \otimes (\nabla'_j)^\vee]) 
 & \rightarrow H^1 (X, \mr{Gr}^j (\V) \otimes (\mcE/\mr{Fil}^j (\V))^\vee) \\
 & \xrightarrow{H^1 (\nabla'''_j \otimes (\nabla'_j)^\vee)} H^1 (X, \Omega_X \otimes (\mr{Gr}^j (\V) \otimes (\V/\mr{Fil}^j (\V))^\vee)). \notag
 \end{align}
 Note that 
 $ \mbH^1 (X, \mcK^\bullet [\nabla'''_j \otimes (\nabla'_j)^\vee])$ (resp., $H^1 (X, \mr{Gr}^j (\V) \otimes (\mcE/\mr{Fil}^j (\V))^\vee)$) classifies extensions $(\mcG, \nabla_\mcG)$ (resp., $\mcG$) of $(\V/\mr{Fil}^j (\V), \nabla'_j)$ (resp., $\V/\mr{Fil}^j (\V)$) by $(\mr{Gr}^j (\V), \nabla'''_j)$ (resp., $\mr{Gr}^j (\V)$) and the first arrow in \eqref{Eq4421}
arises from the assignment between extensions $(\mcG, \nabla_\mcG) \mapsto \mcG$. 
The extension $\varepsilon$ specifies an element $\gamma_\varepsilon$ of $H^1 (X, \mr{Gr}^j (\V) \otimes (\V/\mr{Fil}^j (\V))^\vee)$.
 By the definition of $\nabla''_j$,
 the image of $\gamma_\varepsilon$  via the second arrow  in \eqref{Eq4421} vanishes.
 Hence, there exists an element of $ \mbH^1 (X, \mcK^\bullet [\nabla'''_j \otimes (\nabla'_j)^\vee])$ mapped to $\gamma_\varepsilon$ via the first arrow; this specifies a connection $\nabla'_{j-1}$ on $\V/\mr{Fil}^{j-1}(\V)$ fitting into a short exact sequence
 \begin{align}
 0 \rightarrow (\mr{Gr}^j (\V), \nabla'''_j) \rightarrow (\V/\mr{Fil}^{j-1}(\V), \nabla'_{j-1}) \rightarrow (\V/\mr{Fil}^j (\V), \nabla'_j) \rightarrow 0.
 \end{align}
 This completes the proof of the induction step.
 In particular, by considering the case of  $j=0$,
 we obtain a connection on $\V \left(= \V/\mr{Fil}^0 (\V) \right)$, as desired.
 \end{proof}
\SSP

\begin{rem}[Finite level vs. Infinite level] \label{Rem994}
Considering the structure of the stratified fundamental group $\pi_1^\mr{str}(X)$ (cf. Remark \ref{Reeem441}), we find  that   the only   indecomposable vector bundles admitting  $\mcD^{(\infty)}$-module structures are, up to tensoring with line bundles,  those isomorphic to $\mcF_r$'s
(cf.  ~\cite[Theorem 21]{dSa} and its proof).
This fact is reminiscent of the  characteristic $0$ setting.
 The key  difference, however, lies in the case of finite-level $\mcD$-modules,
 where    a wide variety of vector bundles admit such structures, 
   as suggested by the above proposition.
 This fundamental distinction highlights the necessity of developing specialized arguments for finite-level $\mcD$-modules, setting them apart from their $\mcD^{(\infty)}$-module counterparts.
\end{rem}

\LSP
\subsection{The existence of complete flags} \label{SS44}

In this subsection, we complete the proof of  the remaining part of (a) $\Rightarrow$ (b) in Theorem \ref{ThA}.
To begin with, we establish the following assertion.

\SSP
\ble \label{Lem443}
Let $\V$ and $\V'$ be indecomposable vector bundles on $X$, and write 
$L (\V) = (\mcL_1, \cdots, \mcL_\ell)$, $L (\V')= (\mcL'_1, \cdots, \mcL'_{\ell'})$.
Also, let
 $f : \mcE \rightarrow \mcE'$ be  an $\mcO_X$-linear morphism with $f (\mr{Fil}^1 (\mcE)) \not\subseteq \mr{Fil}^1 (\mcE)$.
 Then,  we have $H^0 (X, {\mcL'_1}^\vee \otimes \mcL_1) = 0$, or equivalently, 
 one of the following two conditions holds:
 \begin{itemize}
 \item
 $\mr{deg}(\mcL_1) < \mr{deg}(\mcL'_1)$;
 \item
$\mr{deg}(\mcL_1) = \mr{deg}({\mcL'_1})$ and $\mcL_1 \not\cong \mcL'_1$.
 \end{itemize}
\ele
\begin{proof}
Suppose, on the contrary, that  $H^0 (X, {\mcL'_1}^\vee \otimes \mcL_1) \neq  0$.
Since $f (\mr{Fil}^1 (\mcE)) \not\subseteq \mr{Fil}^1 (\mcE')$, there exists an injection $\iota : \mcL_1 \hookrightarrow \mr{Fil}^1 (\V)$
such that  $\iota' := \pi \circ f \circ \iota : \mcL_1 \rightarrow \V'/\mr{Fil}^1 (\V')$ is nonzero,
where $\pi$ denotes the natural quotient $\V' \twoheadrightarrow \V'/\mr{Fil}^1 (\V')$.
In particular, we obtain an extension
\begin{align} \label{Eq2244}
0 \rightarrow \mr{Fil}^1 (\V') \rightarrow \pi^{-1} (\mr{Im} (\iota')) \rightarrow \mcL_1 (\cong \mr{Im}(\iota'))\rightarrow 0.
\end{align}
Recall  that
the endomorphism $\mr{ext}_{\varepsilon'}$ of  $H^0 (X, {\mcL'_1}^\vee \otimes (\V'/\mr{Fil}^1 (\mcE')))$
determined, via  \eqref{Eq223},  by 
 the  extension $\varepsilon' : 0 \rightarrow \mr{Fil}^1 (\V') \rightarrow \V' \rightarrow \V'/\mr{Fil}^1 (\V') \rightarrow 0$  under a natural identification 
 \begin{align}
 \mr{Fil}^1 (\V') \cong H^0 (X, {\mcL'_1}^\vee \otimes (\V'/\mr{Fil}^1 (\V'))) \otimes_k \mcL'_1
 \end{align}
  coincides with the identity morphism.
Hence, the composite of  $\mr{ext}_{\varepsilon'}$
with the morphism 
\begin{align}
\left( 0 \neq \right) H^0 (X, {\mcL'_1}^\vee \otimes \mcL_1) \rightarrow H^0 (X, {\mcL'_1}^\vee \otimes (\V'/\mr{Fil}^1 (\V')))
\end{align}
 induced from
 $\iota' : \mcL_1 \rightarrow \V'/\mr{Fil}^1 (\V')$ is nonzero.
Since  \eqref{Eq2244}
corresponds to this composite via 
\begin{align}
\mr{Ext}^1_{\mcO_X} (\mcL_1, \mr{Fil}^1 (\V')) & \xrightarrow{\sim}\mr{Ext}_{\mcO_X}^1(\mcL_1, H^0 (X, {\mcL'_1}^\vee \otimes \V') \otimes_k \mcL'_1)
\\
 & \xrightarrow{\sim} \mr{Ext}_{\mcO_X}^1(\mcL_1, H^0 (X, {\mcL'_1}^\vee \otimes (\V'/\mr{Fil}^1 (\V'))) \otimes_k \mcL'_1)  \\
&  \xrightarrow{\sim}
\mr{Hom}_k (H^0 (X, {\mcL'_1}^\vee \otimes \mcL_1), H^0 (X, {\mcL'_1}^\vee \otimes (\V'/\mr{Fil}^1 (\V'))))
\end{align}
(cf. \eqref{Eq223}, \eqref{Eq3333}),
this extension  turns out to be  nontrivial.
But,  it contradicts the fact that  $f \circ \iota: \mcL_1 \hookrightarrow \pi^{-1} (\mr{Im}(\iota')) \left(\subseteq \mcE' \right)$  specifies a split injection of   \eqref{Eq2244}.
Thus, our initial assumption must be false, completing the proof of this assertion.
\end{proof}
\SSP

The above lemma will be applied in the proof of  the following proposition.

 \SSP
\bpr \label{Lem1}
Let $\V$ and $\V'$ be 
indecomposable vector bundles on $X$, and write 
$L (\V) = (\mcL_1, \cdots, \mcL_\ell)$, $L (\V')= (\mcL'_1, \cdots, \mcL'_{\ell'})$.
Also, let
 $f : \mcE \rightarrow \mcE'$ be  an $\mcO_X$-linear morphism.
\begin{itemize}
\item[(i)]
If any of the following conditions (a)-(c) is satisfied, then $f (\mr{Fil}^1 (\V)) = 0$:
\begin{itemize}
\item[(a)]
$\mr{deg}(\mcL_1) > \mr{deg}(\mcL'_1)$;
\item[(b)]
$m (\V) > 0$, $\mr{deg}(\mcL_1) = \mr{deg}(\mcL_2)$,  and $\mcL_1 \not\cong \mcL'_1$;
\item[(c)]
$m (\V) = m (\V') = 0$,  $\mr{deg}(\mcL_1) = \mr{deg}(\mcL_2)$,  and $\mcL_1 \not\cong \mcL'_1$.
\end{itemize}
\item[(ii)]
If  $\mcL_1 \cong \mcL_2$, then 
 $f (\mr{Fil}^1 (\V)) \subseteq \mr{Fil}^1(\V')$.
\end{itemize}
 \epr
\begin{proof}
We first consider assertion (i).
If the condition (a) is fulfilled (which implies $H^0 (X, {\mcL'_1}^\vee \otimes \mcL_1) \neq 0$), then it follows from Lemma \ref{Lem443} that
$f (\mr{Fil}^1 (\V)) \subseteq \mr{Fil}^1(\V')$.
Since
\begin{align}
\mr{Fil}^1 (\V) \cong H^0 (X, \mcL_1^\vee  \otimes \V) \otimes_k \mcL_1 \ \ \ \text{and} \ \ \   
\mr{Fil}^1 (\V') \cong H^0 (X, {\mcL'}_1^\vee \otimes \V') \otimes_k \mcL'_1,
\end{align}
the restriction of $f$ to  $\mr{Fil}^1 (\V)$ can be regarded as an $\mcO_X$-linear morphism 
\begin{align}
H^0 (X, \mcL_1^\vee \otimes  \V) \otimes_k \mcL_1 \rightarrow H^0 (X, {\mcL'}_1^\vee \otimes \V') \otimes_k \mcL'_1.
\end{align}
However, this must be the zero map due to the assumption $\mr{deg}(\mcL_1) > \mr{deg}(\mcL'_1)$, proving that $f (\mr{Fil}^1 (\V)) = 0$.

Next, suppose that the condition (b) holds.
If $m (\V') > 0$, then $\mcL_1 = \mcL'_1 = \mcO_X (s \sigma)$ with $s := \mr{deg}(\mcL_1) = \mr{deg}(\mcL'_1)$, which contradicts the assumption.
It follows that $m (\V') = 0$, meaning that $\V' \cong \mcL'_1 \otimes  \mcF_{r'}$, where $r' := \mr{rank} (\V')$.
Since each graded piece of the filtration $\mr{Fil}^\bullet (\V)$ is isomorphic to $\mcL'_1 \left(\not\cong \mcL_1 \right)$,
we have $\mr{Hom} (\mcL_1, \V') = 0$, which implies $f (\mr{Fil}^1 (\V)) = 0$, as desired.

We shall  impose the assumption (c).
This implies that $\V$ and $\V'$ can be identified with $\mcL_1 \otimes \mcF_r$ and $\mcL'_1 \otimes \mcF_{r'}$ (for some $r, r' \in \mbZ_{\geq 0}$), respectively.
Since the graded pieces of the filtrations  $\mr{Fil}^\bullet (\V)$ and  $\mr{Fil}^\bullet (\V')$ are  isomorphic to $\mcL_1$ and $\mcL'_1$, respectively, 
the latter two conditions in (c) implies  $f = 0$.
This completes the proof of assertion (i).

Finally, assertion (ii) follows directly from Lemma \ref{Lem443}.
\end{proof}
\SSP

By applying the above assertion, we obtain  the following theorem.

\SSP
\bt \label{Th2}
Let   $\M$ be an element of $\mbZ_{\geq 0} \sqcup \{ \infty \}$
 and 
$(\V, \nabla)$  a  $\mcD^{(\M)}$-bundle of rank $r \in \mbZ_{> 0}$.
Then,
there exists a complete flag on  $(\V, \nabla)$.
\et
\begin{proof}
First, we shall prove the assertion for $m=0$ 
 by induction on $r$.
In the base step, i.e., the case where $r = 1$, there is nothing to prove.
To consider the induction step, we assume that
the assertion for the rank $r -1$ case (where $r  >1$)   has been proved. 
There exists  a collection of line bundles   $\mcN_1, \cdots, \mcN_a$  on $X$ satisfying the following conditions:
\begin{itemize}
\item
$\mcN_1, \cdots, \mcN_a$ are pairwise non-isomorphic;
\item
A line bundle $\mcL$ is isomorphic to $\mcN_{a'}$ for some   $a' \in \{ 1, \cdots, a \}$
if and only if there exists 
an  indecomposable component $\V'$ of $\V$ with $\mcL \cong \mcL_1$, 
where  $L (\V') = (\mcL_1, \cdots, \mcL_\ell)$.
\end{itemize}
After possibly rearranging the order of  indices, 
 one can assume that   $\mr{deg}(\mcN_1) \geq \mr{deg}(\mcN_2) \geq \cdots \geq \mr{deg}(\mcN_a)$ and $\mcN_j \not\cong \mcO_X (d \sigma)$ for $j \neq 1$, where $d := \mr{deg}(\mcN_1)$.
For each $i =1, \cdots, a$,
let $\V_{i, 1}, \cdots, \V_{i, M_i}$ (where $M_i \in \mbZ_{>0}$)   be the  indecomposable components  of $\V$ whose first line bundles in $L (-)$ are isomorphic to $\mcN_1$.
In particular, $\V$ decomposes as  a direct sum 
\begin{align} \label{Eq22}
\V = \bigoplus_{i=1}^\ell \bigoplus_{j =1}^{M_i}  \V_{i, j}^{\oplus e_{i, j}},
\end{align}
where  $e_{i, j} \in \mbZ_{>0}$.
Since $\V$ admits a connection (i.e., $\nabla$), each direct summand $\V_{i, j}$ admits a connection $\nabla_{\V_{i, j}}$ (cf. ~\cite[Lemma 2.4]{Bis1}).
The direct sum
\begin{align}
\nabla' := \bigoplus_{i=1}^\ell \bigoplus_{j =1}^{M_i} \nabla_{\V_{i, j}}
\end{align}
 defines 
a connection on  the right-hand side of  \eqref{Eq22}.
By regarding  $\nabla'$ as a connection on $\V$ via \eqref{Eq22}, 
we obtain 
 the difference $\nabla - \nabla'$,
 which is $\mcO_X$-linear.
 Hence, 
  $\nabla$ can be expressed as 
\begin{align}
\nabla = \nabla' + \sum_{i,  i'} \sum_{j, j'} h_{i, i', j, j'}
\end{align}
for some $h_{i, i', j, j'} \in \mr{Hom}_{\mcO_X} (\V_{i, j}^{\oplus e_{i, j}}, \V_{i', j'}^{\oplus e_{i', j'}}) \left(\subseteq \mr{End}_{\mcO_X} (\V) \right)$ under   the fixed identification $\mcO_X = \Omega_X$.
According to Proposition  \ref{Lem1}, (i) and (ii),
 the sum $\sum_{i,  i'} \sum_{j, j'} h_{i, i', j, j'}$  restricts to an endomorphism of $\bigoplus_{j=1}^{M_1} \mr{Fil}^1 (\mcE_{1, j})$.
Combining this fact with Proposition  \ref{Lem47},
we see that  $\bigoplus_{j=1}^{M_1} \mr{Fil}^1 (\mcE_{1, j})$ is closed under $\nabla$. 
We denote the resulting connection on $\bigoplus_{j=1}^{M_1} \mr{Fil}^1 (\mcE_{1, j})$ by $\nabla''$.
It follows from Proposition \ref{Lem45} that $p \mid \mr{deg}(\mcN_1)$, and that $\mcN_1$ admits a connection $\nabla_{\mcN_1}$ (cf. ~\cite[Theorem 1.1]{Bis1}, ~\cite[Theorem 1.1]{BiSu}).
Hence, after tensoring $(\V, \nabla)$ with the dual of $(\mcN_1, \nabla_{\mcN_1})$,
we may assume that $\mcN_1 = \mcO_X$.
This assumption yields  an identification
\begin{align}
\bigoplus_{j=1}^{M_1} \mr{Fil}^{1} (\V_{1, j}) \cong \mcO_X^{\oplus s}
\end{align}
 for some  $s \in \mbZ_{> 0}$.
Under this identification,
the connection $\nabla''$ can be expressed as $\nabla'' = d + A$ for some 
$s \times s$ matrix $A \in M_s (k)$.
By applying a triangulation of $A$, we find a vector $\vec{v} \in k^{\oplus s} \left(= H^0 (X, \mcO_X^{\oplus s}) \right)$ with  $\nabla'' (\mcO_X \cdot  \vec{v}) \subseteq \Omega_X \otimes (\mcO_X \cdot  \vec{v})$.
The line subbundle $\mcL$ of $\left(\bigoplus_{j} \mr{Fil}^1 (\mcE_{1, j})  \subseteq \right)\mcE$
 determined by $\mcO_X \cdot \vec{v}$ satisfies $\nabla (\mcL) \subseteq \Omega_X \otimes \mcL$, and 
$\overline{\V} := \V/\mcL$ admits a connection $\overline{\nabla}$ induced naturally from $\nabla$.
By the induction hypothesis, $(\overline{\V}, \overline{\nabla})$ admits a complete flag 
$\{ \overline{\V}^j \}_{j=0}^{r-1}$.
Using this filtration, we define  $\V^0 := 0$ and $\V^j := \pi^{-1} (\overline{\V}^{j-1})$ ($j=1, \cdots, r$), where $\pi$ denotes the natural quotient $\V \twoheadrightarrow \overline{\V}$.
The resulting filtration   $\{ \V^j \}_{j=0}^{r}$ specifies a complete flag $(\V, \nabla)$, completing the proof for $\M = 0$.

Also, by Proposition \ref{Prop3200} in the case where the pair ``$(\M, s)$" is taken to be $(0, \M)$,
the assertion for  $\M \in \mbZ_{> 0}$ is  induced from    the assertion just proved.
 
 Finally,  the assertion for $\M = \infty$ follows immediately  from  the structure   of $\pi_1^\mr{str}(X)$ (cf. Remark \ref{Reeem441})  discussed in 
 ~\cite[Theorem 21]{dSa} and its proof.
 Indeed, this result states  that any $\mcD^{(\infty)}$-bundle can be written as a direct sum of the form $\bigoplus_i (\mcL_i, \nabla_{\mcL_i}) \otimes (\V_i, \nabla_{\V_i})$, where
 $(\mcL_i, \nabla_{\mcL_i})$ has  rank $1$ and $(\V_i, \nabla_{\V_i})$ is unipotent, i.e.,   obtained as a successive extension of $(\mcO_X, \nabla_\mr{triv}^{(\infty)})$.
 Since unipotent connections always admit complete flags, the existence of a complete flag for any $\mcD^{(\infty)}$-bundle follows immediately, completing the proof.
\end{proof}

\LSP
\section{The case of hyperbolic  curves} \label{S453}
\LSP

This final section addresses  the remaining part of Theorem \ref{ThA}, specifically the implication (b) $\Rightarrow$ (a).
To construct a flat vector bundle on a hyperbolic curve that does not admit complete flags,
we will leverage a characteristic-$p$ analogue of the Hitchin fibration, known as  the $p$-Hitchin morphism.

\LSP
\subsection{The $p$-Hitchin morphism} \label{SS44r}

Suppose that $g \geq 2$, and 
let  $(\V, \nabla)$ be a flat vector bundle  on $X$ of rank $r \in \mbZ_{> 0}$.
The $p$-curvature $\psi (\nabla)$
of $\nabla$ is regarded as a global section of $F_{X/k}^* (\Omega_{X^{(1)}}) \otimes \mcE nd (\V)$,  where $\mcE nd (\V) := \mcE nd_{\mcO_X} (\V)$.
We associate  to $\psi (\nabla)$ its characteristic polynomial 
\begin{align}
\mr{Char} (\V, \nabla) := \mr{det} (t -\psi (\nabla)) = t^r + a_{r-1} t^{r-1} + \cdots + a_0
\end{align}
with coefficients  $a_i \in H^0 (X, F_{X/k}^{*}(\Omega_{X^{(1)}}^{\otimes (r-i)}))$ ($0 \leq i \leq r-1$).
For each $i$,
the morphism 
\begin{align}
H^0 (X^{(1)}, \Omega_{X^{(1)}}^{\otimes (r-i)}) \rightarrow H^0 (X, F_{X/k}^{*}(\Omega_{X^{(1)}}^{\otimes (r-i)}))
\end{align}
 induced from $\Omega_{X^{(1)}}^{\otimes (r-i)} \rightarrow F_{X/k*} (F_{X/k}^{*}(\Omega_{X^{(1)}}^{\otimes (r-i)}))$ is injective, so  
we regard $H^0 (X^{(1)}, \Omega_{X^{(1)}}^{\otimes (r-i)})$ as a subspace of  $ H^0 (X, F_{X/k}^{*}(\Omega_{X^{(1)}}^{\otimes (r-i)}))$.
Moreover, it is well-known that $a_i$ lies in $H^0 (X^{(1)}, \Omega_{X^{(1)}}^{\otimes (r-i)})$ (cf. ~\cite[Proposition 3.2]{LaPa}).

Denote by $\mcM_{\mr{dR}}$ the moduli stack classifying flat vector bundles on $X$ of rank $r$, and write 
\begin{align}
B^{(1)} := \mr{Spec} \left(\mr{Sym} \left(\bigoplus_{i=1}^{r} H^0 (X^{(1)}, \Omega_{X^{(1)}}^{\otimes i}) \right)^\vee \right).
\end{align}
We identify $B^{(1)}$ with  the moduli space   classifying  polynomials of the form  $t^r + a_{r-1}t^{r-1} + \cdots + a_0$ with $a_i \in H^0 (X^{(1)}, \Omega_{X^{(1)}}^{\otimes (r-i)})$.
The  assignment $(\V, \nabla) \mapsto \mr{Char} (\V, \nabla)$ determines  a $k$-morphism
\begin{align}
\mr{Hitch} : \mcM_\mr{dR} \rightarrow B^{(1)},
\end{align}
which is called the {\bf $p$-Hitchin morphism} (cf.  ~\cite{BeBr}, ~\cite{ChZh1}, ~\cite{ChZh2}, ~\cite{Gro},  ~\cite{JoPa}, ~\cite{LaPa}, ~\cite{Moc}, and  ~\cite{Wak8}). 
The image of this morphism contains  a dense open subscheme of $B^{(1)}$ (cf. ~\cite[Proposition 4.3, Lemma 4.8]{BeBr}).

Next, 
we shall write
\begin{align}
D^{(1)}:= \mr{Spec}\left( \mr{Sym} \left(H^0 (X^{(1)}, \Omega_{X^{(1)}})^{\oplus r}\right)^\vee\right).
\end{align}
This space parametrizes collections of $r$-tuples $(b_1, \cdots, b_r) \in H^0 (X^{(1)}, \Omega_{X^{(1)}})^{\oplus r}$.
The assignment $(b_1, \cdots, b_r) \mapsto \prod_{i=1}^r (t -b_i)$ determines a $k$-morphism
$\gamma : D^{(1)} \rightarrow B^{(1)}$.
This morphism  captures a special subclass of characteristic polynomials corresponding to cases where the $p$-curvature is diagonalizable in a suitable basis.
Observe that 
\begin{align}
\mr{dim}_k(\bigoplus_{i=1}^r H^0 (X^{(1)}, \Omega_{X^{(1)}}^{\otimes i})) &
= \sum_{i=1}^{r} \mr{dim}_k (H^0 (X^{(1)}, \Omega_{X^{(1)}}^{\otimes i})) \\
& = g + \sum_{i=2}^r (2i-1) (g-1) \notag \\
& =  g + (r^2 -1) (g-1) \notag \\
& > rg \notag \\
& = \mr{dim}_k (H^0 (X^{(1)}, \Omega_{X^{(1)}})^{\oplus r}).
\end{align}
Hence, 
 the morphism  $\gamma$ is not dominant,  which   implies  $\mr{Hitch}^{-1} ( B^{(1)} \setminus \mr{Im} (\gamma)) \neq \emptyset$.

\SSP
\bpr \label{Prop931}
Let  $q$ be a $k$-rational point  of $\mr{Hitch}^{-1} ( B^{(1)} \setminus \mr{Im} (\gamma))$, and 
 denote by  $(\V, \nabla)$ the flat vector bundle 
 classified by $q$.
Then, $(\V, \nabla)$ does not admit any complete flag.
\epr
\begin{proof}
Suppose, on the contrary, that there exists a complete flag $\{ \V^j \}_{j=0}^r$  on $(\V, \nabla)$.
For each $i \in \{0, 1\}$,
let
$\mcE nd (\V)^j$
 denote the subsheaf of 
$\mcE nd (\V)$ 
 consisting of endomorphisms 
$h$ of $\V$ with $h (\V^j) \subseteq \V^{j-i}$ for every $i$.
Since each $\V^j$ is preserved by $\nabla$, 
the connection $\nabla^\vee \otimes \nabla$ on $\mcE nd (\V)$ restricts to a connection 
$\nabla^0$ (resp., $\nabla^1$) on $\mcE nd (\V)^0$ (resp., $\mcE nd (\V)^1$).
We obtain a composite of natural morphisms 
\begin{align} \label{eq201}
\mcE nd (\V)^0  \twoheadrightarrow \mcE nd (\V)^0 /\mcE nd (\V)^1 \xrightarrow{\sim} \bigoplus_{j=1}^r \mcE nd (\mcE^j/\mcE^{j-1}) \xrightarrow{\sim}  \mcO_X^{\oplus r},
\end{align}
which  is compatible with the respective connections $\nabla^0$ and 
$(\nabla_\mr{triv}^{(0)})^{\oplus r}$.
The $p$-curvature  $\psi (\nabla)$  belongs to
 $H^0 (X, F_{X/k}^{*} (\Omega_{X^{(1)}})\otimes \mcE nd (\V)^0)$, and 
 its image under 
 the morphism
 \begin{align}
 H^0 (X, F_{X/k}^{*} (\Omega_{X^{(1)}})\otimes \mcE nd (\V)^0) \rightarrow H^0 (X, F_{X/k}^{*} (\Omega_{X^{(1)}}))^{\oplus r}
 \end{align}
induced from  the tensor product of  \eqref{eq201}  and the identity morphism of $F_{X/k}^* (\Omega_{X^{(1)}})$
determines   an $n$-tuple $(b_1, \cdots, b_n)$ of elements in $H^0 (X^{(1)}, \Omega_{X^{(1)}})$ $\left(\subseteq H^0 (X, F^*_{X/k} (\Omega_{X^{(1)}})) \right)$.
It follows that the characteristic polynomial $\mr{Char} (\V, \nabla)$ is given by  $\prod_{i=1}^n (t - b_i)$, i.e., 
coincides with $\gamma (b_1, \cdots, b_n) \in \mr{Im} (\gamma)$.
This contradicts the assumption $\mr{Hitch}^{-1} ( B^{(1)} \setminus \mr{Im} (\gamma))$,  thus completing the proof 
of this proposition.
\end{proof}

\LSP
\subsection{The non-existence of complete flags} \label{SS4f4r}

The following statement completes the remaining part of Theorem \ref{ThA}.

\SSP
\bt \label{Cor32}
Let $m$ be an element of $\mbZ_{\geq 0} \sqcup \{ \infty \}$ and $r$ a positive integer.
Then, there exists a  $\mcD^{(\M)}$-bundle of rank $r$ that does not admit any complete flag. 
\et
\begin{proof}
The assertion for $\M = 0$ is a direct consequence of Proposition \ref{Prop931}, and this assertion  also induces  the case of $\M \in \mbZ_{> 0}$ by applying  Proposition \ref{Prop3200} in the case where the pair ``$(\M, s)$" is taken to be $(0, \M)$.

Next, we shall consider $\M = \infty$.
It follows from ~\cite[Theorem 1.2]{DuMe}
that there exists a stable vector bundle $\V$ on $X$ admitting an isomorphism  $\nu : F_X^{s*}(\V) \xrightarrow{\sim} \V$ for some $s \in \mbZ_{> 0}$.
In particular, since $\nu$ yields  an $F$-divided sheaf structure  on $\V$ (cf. Remark \ref{Reeem441}),
$\V$ carries  a $\mcD^{(\infty)}$-module structure $\nabla$.
Note that $(\V, \nabla)$ is irreducible.
Indeed, if $(\V, \nabla)$ has a proper  $\mcD^{(\infty)}$-subbundle
 $(\V', \nabla')$, then (by considering, e.g.,  the $F$-divided structure  on  its determinant) we see that  $\V'$ must have  degree $0$.
 However, it contradicts the stability of  $\V$.
 In particular, $(\V, \nabla)$ does not admit any complete flag.
 This completes the proof for $\M = \infty$. 
\end{proof}
\SSP

\begin{rem}[Other examples using  dormant opers] \label{Rem4332}
In the proof of the above assertion for $\M = \infty$,
we constructed   an irreducible $\mcD^{(\infty)}$-bundle whose underlying vector bundle is stable.
Here, we  provide  an example of  an {\it unstable} bundle (of rank $2$) that does not admit any complete flag.

To align with the notation  in ~\cite{Wak12}, we set $\N := \M +1$.
Fix a theta characteristic $\varTheta$ of $X$, i.e., a line bundle $\varTheta$ equipped with an isomorphism $\eta : \varTheta^{\otimes 2} \xrightarrow{\sim}\Omega_X$.
We shall regard   $\nabla_\mr{triv}^{(\N -1)}$ as a $\mcD^{(\N -1)}$-module structure  on $\mcT_X \otimes \varTheta^{\otimes 2}$ via the identification with $\mcO_X$ arising from $\eta$.
Thus, the pair $\vartheta := (\varTheta, \nabla_\mr{triv}^{(\N -1)})$ defines  a $2^{(\N)}$-theta characteristic, in the sense of
~\cite[Definition 5.12]{Wak12}.
We set $\V_\varTheta := \mcD_{\leq 1}^{(\N -1)} \otimes \varTheta$ and $\V_\varTheta^1 := \mcD_{\leq 0}^{(\N -1)} \otimes \varTheta \left(\subseteq \V_\varTheta \right)$.

 Following ~\cite[Definition 5.15]{Wak12}, 
a {\bf $(\mr{GL}_2^{(\N)}, \vartheta)$-oper} on $X$ is defined as   a $\mcD^{(\N -1)}$-module structure $\nabla$ on $\V_\varTheta$ satisfying the following two conditions:
\begin{itemize}
\item
The $\mcO_X$-linear morphism $\mcD^{(\N -1)} \otimes \V_\varTheta \rightarrow \V_\varTheta$ given by $\nabla$ restricts to an isomorphism $\mcD_{\leq 1}^{(\N -1)} \otimes \V_\varTheta^1 \xrightarrow{\sim} \V_\varTheta$.
\item
The connection $\mr{deg}(\nabla)$ on $\mr{det}(\V_\varTheta)$ induced naturally  from $\nabla$ coincides with $\nabla_\mr{triv}^{(\N -1)}$ via the natural composite
\begin{align}
\mr{det}(\V_\varTheta) \xrightarrow{\sim} 
 (\V_\varTheta/\V_\varTheta^1) \otimes \V_\varTheta \xrightarrow{\sim} (\mcT_X \otimes \varTheta) \otimes \varTheta
\xrightarrow{\sim} \mcT_X \otimes \varTheta^{\otimes 2}.
\end{align}
\end{itemize}
Moreover, 
 a $(\mr{GL}_2^{(\N)}, \vartheta)$-oper $\nabla$ is said to be  {\bf dormant} if it has vanishing $p^\N$-curvature.
For the study of opers in positive characteristic,  we refer the reader to  ~\cite{JoPa}, ~\cite{JRXY},  ~\cite{Moc}, ~\cite{Wak8}. 

The notion of an isomorphism between  $(\mr{GL}_2^{(\N)}, \vartheta)$-opers are defined in a suitable manner (though  we here omit the details).
In particular, we obtain the set $\mcO p_{\N}^{^\mr{Zzz...}}$
consisting of isomorphism classes of
dormant $(\mr{GL}_2^{(\N)}, \vartheta)$-opers on $X$.
Under the assumption $p >2$, we know that 
  $\mcO p_\N^{^\mr{Zzz...}}$ is a nonempty finite set for every $\N \in \mbZ_{> 0}$ (cf. ~\cite[Eq.\,(5.8.2), Theorem 6.16, and Proposition 6.26]{Wak12}),
and  that the limit of the projective system
\begin{align} \label{eq10}
\cdots \rightarrow \mcO p_\N^{^\mr{Zzz...}} \rightarrow \cdots \rightarrow \mcO p_2^{^\mr{Zzz...}} \rightarrow \mcO p_1^{^\mr{Zzz...}}
\end{align}
obtained by reducing levels is  nonempty (cf. ~\cite[Corollary 8.28]{Wak12}).
By choosing an element of this limit,
we obtain a $\mcD^{(\infty)} \left(= \varinjlim_{\N \in \mbZ_{>0}} \mcD^{(\N)}\right)$-module structure on $\V_\varTheta$.

Note that $(\V_\varTheta, \nabla)$ is irreducible.
Indeed, suppose, on the contrary, that it admits a line $\mcD^{(\infty)}$-subbundle $(\mcL, \nabla_\mcL)$.
Then, $\mcL$ must have   degree $0$, so
the natural  composite $\mcL \hookrightarrow \V_\varTheta \rightarrow \V_\varTheta / \V_\varTheta^1 \left(\cong \varTheta^{\vee} \right)$ coincides with  the zero map because $0 = \mr{deg}(\mcL) > \mr{deg} (\varTheta^{\vee}) = 1-g$.
This means that $\mcL$ is contained in $\V_\varTheta^1$.
Since $0 = \mr{deg}(\mcL) < \mr{deg}(\V_\varTheta^1) = g-1$,
this contradicts the assumption that $\mcL$ is a subbundle of $\V_\varTheta$.
Consequently, $(\V_\varTheta, \nabla)$ is irreducible, and  in particular,  does not admit any complete flag.
 \end{rem}

\LSP
\subsection*{Acknowledgements}
The second author was partially supported by 
 JSPS KAKENHI Grant Number 21K13770.

\end{document}